\documentclass[12pt]{amsart}
\usepackage{amsmath}
\usepackage{amssymb}
\usepackage{amsfonts}
\usepackage{amssymb,amsthm,amsmath,hyperref}
\usepackage{tikz}
\usepackage{graphicx,float}

\numberwithin{equation}{section}

\newtheorem{example}{Example}[section]
\newtheorem{thm}{Theorem}[section]

\theoremstyle{definition}
\newtheorem{defn}{Definition}[section]
\theoremstyle{remark}
\newtheorem{rem}{Remark}[section]

\newcommand{\ndF }{\mathbb{F}}

\allowdisplaybreaks

\newcommand{\Dtriangle}[7]{
\rule[-3\unitlength]{0pt}{12\unitlength}
\begin{picture}(18,7)(0,3)
\put(4,4){\ifthenelse{\equal{#1}{l}}{\circle*{2}}{\circle{2}}}
\put(5,4){\line(1,0){8}}
\put(14,4){\ifthenelse{\equal{#1}{r}}{\circle*{2}}{\circle{2}}}
\put(4.4472,4.8944){\line(1,2){4.1056}}
\put(9,14){\ifthenelse{\equal{#1}{t}}{\circle*{2}}{\circle{2}}}
\put(13.5528,4.8944){\line(-1,2){4.1056}}
\put(2,3){\makebox[0pt][r]{\scriptsize #2}}
\put(9,16){\makebox[0pt]{\scriptsize #3}}
\put(16,3){\makebox[0pt][l]{\scriptsize #4}}
\put(6,9){\makebox[0pt][r]{\scriptsize #5}}
\put(12.5,9){\makebox[0pt][l]{\scriptsize #6}}
\put(9,1){\makebox[0pt]{\scriptsize #7}}
\end{picture}}

\title{Left symmetric algebras from DNA insertion }

\author{Chen Yuan}
\address{Chen Yuan, Beijing Forestry University, Beijing, 100083, China}
\email{yuanc@bjfu.edu.cn}

\author{Zhixiang Wu}
\address{Zhixiang Wu, Department of Mathematics,
Zhejiang University,
Hangzhou, 310027, China}
\email{wzx@zju.edu.cn}

\author{Jing Wang}\thanks{The work is supported by the National Natural Science Foundation of China (No.12471038)}
\address{Jing Wang,
Beijing Forestry University, Beijing, 100083, China}
\email{wang\_jing619@163.com}

\date{}

\begin{document}
\maketitle
\begin{abstract}
DNA recombination is a fundamental biological process that encodes genetic information for organism development and function. In this study, we construct the left symmetric algebras arising
from the operation of DNA insertion. We define a new operation of
insertion by modifying the simplified insertion
$$x\Rightarrow y:=f(\mid x\mid,\ \mid y \mid)\sum\limits_{i=0}^{q}
y_{1}y_{2}\cdots y_{i} x y_{i+1}\cdots y_{q},$$
where $x =
x_{1}x_{2}\cdots x_{p}$, $y = y_{1}y_{2}\cdots y_{q}$, and $\mid x\mid, \mid y\mid$ denote the lengths of $x$ and $y$, respectively. We prove that the algebra $\mathbb{F}(R)$ (over a field $\mathbb{F}$ of characteristic $0$, with $R$ being an infinite free semigroup generated by DNA nucleotides $\{A, G, C, T\}$) forms a left symmetric algebra if and only if the function $f$ satisfies the condition
$$f(m, n) f(m+n, p)=f(n, p) f(m, n+p)=
f(m, p) f(n, m+p),$$ where $m, n, p\in \mathbb{N}$. A key example of such a function is $f(m, n)=\exp\{g(m, n)\}$, where $g(m,
n)=k\cdot mn,$ and $k$ is a fixed positive number, which effectively models length-dependent DNA insertion dynamics. This work enriches the theory of non-associative algebras and provides a mathematical framework for quantitative analysis of DNA recombination processes.

{\bf{Key Words}}: Left symmetric algebras,  DNA insertion operations, Mathematical Model
\end{abstract}
\section{Introduction}
Genetic information is obtained through DNA recombination. Recombination
 provides the genetic program for the development and
functioning of all living organisms. The algebraic formalization of
DNA recombination is represented in the form of a linear space $\ndF(R) $
over a field $\ndF$ of characteristic 0, where $R$ is an infinite free
semigroup generated by the set of DNA nucleotides $\{A,G,C,T\}$. The
schematic model of non-homologous DNA recombination can be
represented in the form: $c \cdot (ab) \rightarrow a c b $, where
two chromosomes $(a b)$ and $c$ participate in non-homologous
recombination. The algebraic formalization of the above
recombination, with respect to all possible insertions of DNA $c$ in
DNA $(a b)$, defines the algebra of simplified insertions. Let $X\in
R$. Then $X$ can be presented as $X=x_1x_2\cdots x_n$, where
$x_i\in\{A,G,C,T\}$. The operation, defined in \cite{sverchkov2010algebraic}, is given by
\begin{eqnarray}\label{VDD-E1.1}a\cdot X=\sum_{i=0}^nx_1\cdots x_{i}ax_{i+1}\cdots
x_n, \end{eqnarray} where $a$ is an arbitrary element in $\ndF(R)$; and
$a \cdot b =\sum \alpha_s(a\cdot u_s)$ , where $ b =\sum \alpha_su_s$,
$\alpha_s \in \ndF$ , and $u_s$ are monomials in $\ndF(R)$. This operation
$\cdot$ is called the operation of left simplified insertion. And it satisfies the left symmetric identity $(x, y, z) = ( y, x, z) $. (Notice that the categories of the left symmetric and right-symmetric algebras are equivalent \cite{sverchkov2010algebraic}.) The
operation of simplified insertion was first introduced by 
Bremner \cite{bremner2005dna}, and it is an algebraic formalization of the operation
of normal insertion in the theory of DNA computing.
 This algebra defined on
$\ndF(R)$ is a very important class of non-associative algebras. Non-associative algebras
have previously been applied to population genetics in the well-developed theory of genetic algebras. For surveys of this area,
see \cite{reed1997algebraic}.

 Left symmetric algebras (also called pre-Lie algebras) have gained extensive attention and development in recent years and play a significant role in many fields. In 1890, A. Cayley first introduced rooted tree algebras into the study of left symmetric algebras, see \cite {cayley1890theory}. By the 1960s, left symmetric algebras began to emerge in geometric and algebraic contexts such as convex homogeneous cones \cite {vinberg1963convex} and affinely flat manifolds \cite{koszul1961domaines,matsushima1968affine} and the deformation theory of associative algebras \cite{gerstenhaber1963cohomology}. The scholars then began to study left symmetric algebras, focusing on issues such as their structures, classification, and important applications in different fields \cite{burde2006left,WOS:000461394300021,WOS:000781005400004}. In this paper, we construct new left symmetric algebras through DNA recombination to enrich the research field.

The left simplified insertion left symmetric algebra is
infinite-dimensional as a vector space. However, one finds that there are only finitely many types of
DNA in nature. Thus we  modify the operation of the left simplified insertion in this paper.   

The paper is organized as follows. In Section 2, we recall left symmetric algebras and the simplified insertion, and prove that the simplified insertion satisfies the left symmetric identity. Then we show that the intermediate version of the insertion
operation does not
satisfy the left symmetric identity, which answers an open problem in \cite{bremner2005dna}. And we prove that the operation of
synchronized insertion given in \cite{daley2004families} does not
satisfy the left symmetric identity, either. In Section 3, we specify the conditions under which the simplified insertion operation satisfies the left symmetric identity. In Section 4, we
give the definition of modified simplified insertion and construct a new left symmetric structure for the modified simplified insertion. Specifically, We define a new operation of insertion as follows
$$x\Rightarrow y:=f(\mid x\mid, \mid y
\mid)\sum\limits_{i=0}^{q} y_{1}y_{2}\cdots y_{i} x y_{i+1}\cdots
y_{q},$$ where $x = x_{1}x_{2}\cdots x_{p}$, $y =
y_{1}y_{2}\cdots y_{q}$, and $\mid x\mid$, $\mid y\mid$ denote the lengths of $x$ and $y$, respectively. This operation is called modified simplified insertion. Then $\ndF(R)$ is a left symmetric
algebra if and only if $f(m,n)$ satisfies 
$$f(m, n) f(m+n, p)=f(n, p) f(m, n+p)= f(m, p) f(n, m+p),$$ where $m, n, p\in
\mathbb{N}$.
For example, if the operation
is given by 
\begin{eqnarray*}
f(m,n):=e^k(m*n),
\end{eqnarray*}
where $k$ is a fixed negative integer, then $\ndF(R)$ is still a left symmetric algebra. And the
new genetics generated by the left simplified insertion can be
suppressed by this coefficient $e^k(m*n)$ if $m$ or $n$ is large enough. 
At last, we illustrate properties of these left symmetric algebras and certain applications in genetics and molecular genetics.

\section{Definitions and Notions}
\subsection{The left symmetric algebra}

   Let $\ndF$ be a field of characteristic zero. Let $A$ be an algebra over $\ndF$. For elements $x$
   and $y$ in $A$, the bilinear product is denoted by $(x, y) \rightarrow  x\circ y$. Given three elements $x$, $y$ and $z$, we denote the associator of these elements by  $$(x, y, z) = (x\circ y)\circ z - x\circ(y\circ z ).$$
\begin{defn}
  A left symmetric algebra $A$ is an algebra whose associator satisfies
  $$
  \forall x, y, z \in  A,(x, y, z) =(y, x, z)\ \ \text{or equivalently} \ \  (x\circ y)\circ z - x\circ(y\circ z )=(y\circ x)\circ z - y\circ(x\circ z ).
  $$ And the identity $(x, y, z) =(y, x, z)$ is called the left symmetric identity.
\end{defn}

\begin{example}
   \item[(a)] Every associative algebra is a left symmetric algebra with the left symmetric structure given by
   $$x\circ y:= x y.$$ Obviously, $(x\circ y)\circ z - x\circ(y\circ z )=0$.
   \item[(b)] Let $A$ be the vector space $C^{\infty}(\mathbb{R},\mathbb{R})$ of smooth functions. For $f$ and $g$ in $A$, we set $f \circ g := f \frac{dg}{dx}$. It is easy to check that $A$ is a left symmetric algebra.
\end{example}

\subsection{Languages and Simplified insertion}
  Let $S$ be a finite non-empty set, denoted by $S=\{a_i\mid i\in I \}\ $. A word over $S$ is a finite string $w = a_{1}a_{2}\cdot\cdot\cdot a_{p}$, where $p\geq  0$ and $a_{i} \in S (1 \leq i \leq p)$. For $p = 0$
  we have the empty word denoted 1. We denote the length of $w$ by $\mid w\mid = p$. We write $M(S)$ for the set of all words and any subset of $M(S)$ is called a language. 
  
Consider the operation of DNA insertion of one word into another. This gives a method for the combination of DNA molecules that are well studied in molecular genetics.
  Given $x, y \in M(S)$, we consider all insertions of $x$ into $y$:$$x\rightarrow y = \{ y^{1}xy^{2}\mid y = y^{1} y^{2}, y^{1}, y^{2}\in M(S) \}.$$
  Note that we allow both $y^{1}$ and $y^{2}$ to be empty.
  
   Firstly, we recall the simplified insertion.
\begin{defn}
 The simplified insertion of $x$ into $y$ is the linearized form of this operation, that is, the sum of all insertions, as
 $$x\rightarrow y =\sum\limits_{i=0}^{q} y_{1}y_{2}\cdots y_{i} x y_{i+1}\cdots y_{q},$$
 where $y = y_{1}y_{2}\cdots y_{q}$, $y_{i} \in S $ and  $x = x_{1}x_{2}\cdots x_{p}$, $x_{i} \in S$.
\end{defn}

 Let $A[S]$ be the free associative algebra over $\ndF$ generated by $S$. We can define a new operation of multiplication $\circ $ by 
\begin{eqnarray}\label{2VDD-E1.1}
  x\circ y=x\rightarrow y =\sum\limits_{i=0}^{q} y_{1}y_{2}\cdots  y_{i} x y_{i+1}\cdots y_{q},
\end{eqnarray}
  where $y =y_{1}y_{2}\cdots y_{q}$, $y_{i} \in S $. The associator for this new operation, that is simplified insertion operation, is defined as usual by
\begin{eqnarray}\label{2VDD-E1.2}
(x, y, z)_1 = (x\circ y)\circ z - x\circ(y\circ z )=((x\rightarrow y)\rightarrow z) - (x\rightarrow(y\rightarrow z )).
 \end{eqnarray}

\begin{example}
   Let $x =a_{1}a_{2}a_{3},\ y=a_{4}a_{5}$ and $z=a_6$. Then $x\rightarrow y =a_{1}a_{2}a_{3}a_{4}a_{5}+a_{4}a_{1}a_{2}a_{3}a_{5}+a_{4}a_{5}a_{1}a_{2}a_{3}$ and $(x, y, z)_1=-a_1a_2a_3a_6a_4a_5-a_4a_5a_6a_1a_2a_3.$
\end{example}

\begin{thm}\label{2VDD-T1.1}
 The algebra $A[S]$ is a free left symmetric algebra whose associator is defined by $(\ref{2VDD-E1.2})$.
\end{thm}

\begin{proof}
Let $$x = x_{1}x_{2}\cdots x_{p},\  \ y = y_{1}y_{2}\cdots y_{q},\  \ z = z_{1}z_{2}\cdots z_{r}.$$\\
Using $(\ref{2VDD-E1.1})$, we obtain
           $$x\rightarrow y=\sum\limits_{j=0}^{q} y_{1}y_{2} \cdots y_{j} x y_{j+1}\cdots y_{q},$$
           $$y\rightarrow z=\sum\limits_{k=0}^{r} z_{1}z_{2}\cdots z_{k} y z_{k+1}\cdots z_{r}.$$
Then
   \begin{eqnarray}
  &&x\rightarrow(y\rightarrow z)
  =\sum\limits_{k=0}^{r}(\sum\limits_{j=0}^{k-1} z_{1}z_{2}\cdots z_{j}x z_{j+1} \cdots z_{k} y z_{k+1}\cdots z_{r}\nonumber\\&&+ z_{1}z_{2}\cdots z_{k}(x\rightarrow y)z_{k+1}\cdots z_{r}+\sum\limits_{j=k+1}^{r} z_{1}z_{2}\cdots z_{k} y z_{k+1} \cdots z_{j}xz_{j+1}\cdots z_{r})\nonumber,
    \end{eqnarray}
  \begin{eqnarray}
   (x\rightarrow y)\rightarrow z=\sum\limits_{j=0}^{q}\sum\limits_{k=0}^{r} z_{1}z_{2}\cdots z_{k}y_{1}y_{2}\cdots  y_{j} x y_{j+1}\cdots y_{q}z_{k+1}\cdots z_{r}\nonumber.
   \end{eqnarray}
Hence, we obtain
\begin{eqnarray}
           &&((x\rightarrow y)\rightarrow z) - (x\rightarrow(y\rightarrow z))\nonumber\\
           &=&-\sum\limits_{k=0}^{r}(\sum\limits_{j=0}^{k-1} z_{1}z_{2}\cdots z_{j}x z_{j+1} \cdots z_{k} y z_{k+1}\cdots z_{r}\nonumber\\
           &&+\sum\limits_{j=k+1}^{r} z_{1}z_{2}\cdots z_{k} y z_{k+1} \cdots z_{j}xz_{j+1}\cdots z_{r})\nonumber\\&:=& LHS,\label{2VDD-E1.8}
\end{eqnarray}
\begin{eqnarray}
           &&((y\rightarrow x)\rightarrow z) - (y\rightarrow(x\rightarrow z))\nonumber
           \\&=&-\sum\limits_{k=0}^{r}(\sum\limits_{j=0}^{k-1} z_{1}z_{2}\cdots z_{j}y z_{j+1} \cdots z_{k} x z_{k+1}\cdots z_{r}\nonumber\\
           &&+\sum\limits_{j=k+1}^{r} z_{1}z_{2}\cdots z_{k} x z_{k+1} \cdots z_{j}y z_{j+1}\cdots z_{r})\nonumber
           \\&:=& RHS.\label{2VDD-E1.9}
\end{eqnarray}
The terms of $LHS$ are those in which both $x$ and $y$ are inserted into
$z$ and separated by at least one letter of $z$. Hence, $LHS=RHS$. 
We immediately obtain \begin{eqnarray}
  ((x\rightarrow y)\rightarrow z) - (x\rightarrow(y\rightarrow z ))=((y\rightarrow x)\rightarrow z )- (y\rightarrow(x\rightarrow z )).\nonumber
  \end{eqnarray}
\end{proof}

\begin{rem}
  We consider the following open problem given in \cite{bremner2005dna} by Bremner.\\
  Problem: An intermediate version of the insertion operation is obtained when we regard the set-theoretic definition as producing a set rather than a multiset. That is, on words $x = x_1\cdots x_p$ and $y = y_1 \cdots y_q$ with $x_i, y_j \in S,$ we define
  \begin{eqnarray}\delta(x, y)=
\begin{cases}
1, &\text{if}\ x_1 = y_1\cr 0, &\text{otherwise} \end{cases} \nonumber
\end{eqnarray}
We then consider the operation
$$x\circ y=\sum\limits_{j=0}^{q}\delta(x, y_j\cdots y_q) y_{1}y_{2}\cdots y_{j} x y_{j+1}\cdots y_{q}.$$
What are the polynomial identities satisfied by this operation?

We can check that this operation does not satisfy the left symmetric identity $(x, y, z) = (y, x, z)$.
We present the following counterexample. Let $S=\{a, b, c\}, \ x=ab, \ y=abc$ and $ z=ac.$ Then $x\circ y= ababc, \ (x\circ y)\circ z=ababcac, \ x\circ (y\circ z)=ababcac + abcabac$ and hence $(x\circ y)\circ z- x\circ (y\circ z)=- abcabac.$ Similarly, $(y\circ x)\circ z- y\circ (x\circ z)=- ababcac.$ That is $(x, y, z) \not= (y, x, z)$.
\end{rem}

\begin{rem}
  Recall another DNA insertion operation, the linearized version of the synchronized insertion introduced in \cite{daley2004families}. We show that it does not satisfy the left symmetric identity $(x, y, z) = (y, x, z)$, either.
  $$x\rightrightarrows y=\sum\limits_{j=1}^{q}t(x, y_j\cdots y_q) y_{1}y_{2}\cdots y_{j} x y_{j+1}\cdots y_{q},$$
  where $t(x, y)$  is defined by $$t(x, y) = k  \ \ \text{when}\ \  x_i = y_i \ \ \text{for} \ \ 1 \leqslant i \leqslant k, \ \text{but}\ \ x_{k+1} \ne y_{k+1} (\text{or} \ k + 1 \geqslant min(p, q)).$$
  Without loss of generality, we choose the simplest case of the linearized version of synchronized insertion in which $S$ contains one letter: $S=\{a\}.$ Then the synchronized insertion becomes
 $$a^p\rightrightarrows a^q=c(p, q)a^{p+q},$$ where
 \begin{eqnarray} c(p, q)=
 \begin{cases}
 \frac{1}{2}p(2q-p+1), &\text{if}\ p<q\cr \frac{1}{2}q(q+1),  &\text{otherwise} 
 \end{cases} \nonumber
 \end{eqnarray}
 Let $x=a^p, \ y=a^q$ and $z=a^r,$ where $p=2, \ q=3$ and $r=6$. Then $(x\circ y )\circ z- x\circ (y\circ z)-((y\circ x)\circ z- y\circ (x\circ z))=a^{16}\neq 0.$ Hence, it  does not satisfy the left symmetric identity.
 \end{rem}
\section{One of the subspaces and left symmetric algebra}
  Firstly, we describe how signs may be introduced in the notions of words and DNA insertions. Set $S=\{a_{1},a_{2},\cdots,a_{p}\}$. For any $a_{i} , a_{j} \in {S}$, define
  \begin{eqnarray}\delta'(a_{i},a_{j})=
\begin{cases}
1, &\text{if}\ a_{i}\ \text{and}\ a_{j}\  \text{can be adjacent in a word} \cr 0, &\text{otherwise} \end{cases} \nonumber
\end{eqnarray}
 That is, we say that $a_i$ and $a_j$ can be adjacent in a word if there are some words like $\cdots a_ia_j \cdots $ and therefore $\delta'(a_{i},a_{j})=1$; otherwise, $\delta'(a_{i},a_{j})=0$. In particular, we set $\delta' (a_{i},a_{i})= 1, \ \text{for}\ \forall\ a_{i} \in {S} $. 

We give examples to illustrate the role of $\delta'(a_{i},a_{j})$.

\begin{example}
\item[(a)] Let $S=\{a, b, c\}$ and $\delta' (a,b)=\delta' (b,c)=\delta' (a,c)=1$. We can use the adjacency graph \\

\begin{picture}(8,8)
\setlength{\unitlength}{1.5mm}
\put(0,0){\circle{1}}
\put(0.5,0){\line(1,0){6}}
\put(7,0){\circle{1}}
\put(3.6,4){\circle{1}}
\put(0.2,0.3){\line(1,1){3.2}}
\put(3.8,3.5){\line(1,-1){3.2}}
\put(-1.2,-1){\makebox[0pt]{\scriptsize $a$}}
\put(8.2,-1){\makebox[0pt]{\scriptsize $c$}}
\put(4.2,4.6){\makebox[0pt]{\scriptsize $b$}}
\end{picture}  \ \ \ \ \
\vspace{1.5mm}
to illustrate the relationships. Set $x=abc$ and $y=bc$. Then $x\circ y=abcbc+ babcc+ bcabc.$
\item[(b)] Let $S=\{a, b, c\}$, $\delta' (a,b)=\delta' (b,c)=1$ and $\delta' (a,c)=0.$ The adjacency graph is \\
\begin{picture}(8,8)
\setlength{\unitlength}{1.5mm}
\put(0,0.7){\circle{1}}
\put(0.5,0.7){\line(1,0){3.9}}
\put(5,0.7){\circle{1}}
\put(10,0.7){\circle{1}}
\put(5.5,0.7){\line(1,0){3.9}}
\put(-0.8,-1){\makebox[0pt]{\scriptsize $a$}}
\put(5,-1.3){\makebox[0pt]{\scriptsize $b$}}
\put(10.6,-1){\makebox[0pt]{\scriptsize $c$}}
\end{picture} \ \ \ \ \ \ \ \ \
\vspace{1mm}
. Set $x=abc$ and $y=bc$. Then $x\circ y=abcbc+ babcc$.
\end{example}

Notice that $A[S]$ is the associative algebra over $\ndF$ generated by $S$. 
Consider the subspace of $A[S]$ for which there exist $ a_{i}, a_{j} \in S$ such that $\delta' (a_{i},a_{j})= 0$. 

\begin{thm}\label{Th3.1}
   Let $S=\{a_{1}, \ \dots,a_{p}\}$, $ \mid S\mid = p.$
\item[(a)]
   If $\mid S \mid= 1$, then the operation $\circ$ satisfies the left symmetric identity.
\item[(b)]
   If $\mid S\mid = 2$ and $\delta' (a_{1},a_{2})=1 $, then the operation $\circ$ satisfies the left symmetric identity.
   If $\mid S\mid = 2$ and $\delta'(a_{1},a_{2})=0$, then the operation $\circ$ is associative. It obviously satisfies the left symmetric identity.
\item[(c)]
   If $\mid S\mid \geqslant 3 $ and there exist $a_{i}, a_{j} \in S$ such that $\delta'(a_{i},a_{j})=0 $, then the operation $\circ$ is non-associative and non-commutative. It does not satisfy the left symmetric identity.
   Otherwise, the operation is non-associative and non-commutative. However, it satisfies the left identity.
\end{thm}

\begin{proof}
  If $\delta' (a_{i},a_{j})= 1$ for all $ a_{i}, a_{j} \in S$ then the operation satisfies the left symmetric identity by Theorem \ref{2VDD-T1.1}. 
  If $\mid S \mid = 2$ and $\delta'(a_{1},a_{2})=0 $, then the operation obviously satisfies $(a_{i}\circ a_{j})\circ a_{k}= a_{i} \circ(a_{j}\circ a_{k})= 0 $ for $a_{i}, a_{j}, a_{k} \in \{a_{1},a_{2}\}$. Hence, we have $(a)$ and $(b)$. 
  It is not difficult to prove that the claims are true for $p = 3$.
  $(c)$ follows since the space of $p=3$ is the subspace of $p \geqslant 3$.
\end{proof}

\begin{example}
 Let $S=\{a, b, c\}$, $\delta' (a,b)=\delta' (b,c)=1$ and $\delta' (a,c)=0$. Set $x=a, \ y=b$ and $z=c$, then $x\circ y=ab+ ba $ and $y\circ z= bc+ cb $.\\
 Then
\begin{eqnarray}
           &&(x\circ y)\circ z =abc+ cba, \ (y\circ x)\circ z = abc+ cba,\nonumber\\
           &&x\circ (y\circ z) = abc+ cba, \ y\circ (x\circ z) =0,\nonumber\\
           &&LHS:=(x\circ y)\circ z-x\circ (y\circ z)=0,\nonumber\\
           &&RHS:=(y\circ x)\circ z-y\circ (x\circ z)=abc+ cba.\nonumber
\end{eqnarray}
Hence, it does not satisfy the left symmetric identity.
\end{example}


\section{The modified simplified insertion}
 \subsection{The modified simplified insertion}
 \begin{defn}
 We define the modified simplified insertion as follows
 $$x\Rightarrow y:=f(\mid x\mid, \mid y \mid)\sum\limits_{i=0}^{q} y_{1}y_{2}\cdots y_{i} x y_{i+1}\cdots y_{q},$$
 where the function $f: \mathbb{N}\times\mathbb{N}\longrightarrow \mathbb{R}$ is given by $(m, n)\mapsto f(m, n)$.
 Then we have
\begin{eqnarray}\label{2VDD-E1.3}
(x, y, z)_2 = (x\circ y)\circ z - x\circ(y\circ z )=((x\Rightarrow y)\Rightarrow z) - (x\Rightarrow(y\Rightarrow z )).
 \end{eqnarray}
\end{defn}

 \begin{example}
    Let $x =a_{1}a_{2}a_{3}$ and $y=a_{4}a_{5}$, then $x\Rightarrow y =f(3,2)(a_{1}a_{2}a_{3}a_{4}a_{5}+a_{4}a_{1}a_{2}a_{3}a_{5}+a_{4}a_{5}a_{1}a_{2}a_{3})$.
\end{example}

\subsection{The left symmetric algebra for modified simplified insertion}

 \begin{thm}\label{2VDD-T1.2}
  $A[S]$ is the left symmetric algebra where $(x, y, z)$ is defined as $(\ref{2VDD-E1.3})$ if and only if $f$ satisfies the following conditions:
  \begin{eqnarray}\label{2VDD-E2.0}
  f(m, n) f(m+n, p)=f(n, p) f(m, n+p)= f(m, p) f(n, m+p),
  \end{eqnarray}
  where $m, n, p\in \mathbb{N}$.
\end{thm}

\begin{proof}
  Let $$x = x_{1}x_{2}\cdots x_{p},\  \ y = y_{1}y_{2}\cdots y_{q},\  \ z = z_{1}z_{2}\cdots z_{r}.$$\\
 We obtain
           $$x\Rightarrow y=f(\mid x \mid, \mid y \mid)\sum\limits_{j=0}^{q} y_{1}y_{2} \cdots y_{j} x y_{j+1}\cdots y_{q},$$
           $$y\Rightarrow z=f(\mid y \mid, \mid z \mid)\sum\limits_{k=0}^{r} z_{1}z_{2}\cdots z_{k} y z_{k+1}\cdots z_{r}.$$
  Then
   \begin{eqnarray}
  &&x\Rightarrow(y\Rightarrow z)\nonumber\\
  &=&f(\mid y \mid,\mid z \mid)f(\mid x \mid,\mid y\mid +\mid z \mid)\sum\limits_{k=0}^{r}(\sum\limits_{j=0}^{k-1} z_{1}z_{2}\cdots z_{j}x z_{j+1} \cdots z_{k} y z_{k+1}\cdots z_{r}\nonumber\\
  &&+ z_{1}z_{2}\cdots z_{k}(x\rightarrow y)z_{k+1}\cdots z_{r}+\sum\limits_{j=k+1}^{r} z_{1}z_{2}\cdots z_{k} y z_{k+1} \cdots z_{j}x z_{j+1}\cdots z_{r})\nonumber,\label{2VDD-E2.1}
    \end{eqnarray}
  \begin{eqnarray}
   &&(x\Rightarrow y)\Rightarrow z\nonumber\\
   &=&f(\mid x \mid,\mid y \mid)f(\mid x\mid+\mid y \mid,\mid z \mid)\sum\limits_{j=0}^{q}\sum\limits_{k=0}^{r} z_{1}\cdots z_{k}y_{1}\cdots  y_{j} x y_{j+1}\cdots y_{q}z_{k+1}\cdots z_{r}\nonumber.\label{2VDD-E2.2}
   \end{eqnarray}
Denote $H(\mid x\mid,\mid y\mid,\mid z\mid)$ and $H_2(\mid x\mid,\mid y\mid,\mid z\mid)$ as $$f(\mid x \mid,\mid y \mid)f(\mid x\mid +\mid y \mid,\mid z \mid)-f(\mid y \mid,\mid z \mid) f(\mid x \mid,\mid y\mid+\mid z \mid))$$ and $$f(\mid y \mid,\mid z\mid)f(\mid x \mid,\mid y\mid +\mid z \mid),$$respectively.
Hence, we obtain
\begin{eqnarray}
           &&((x\Rightarrow y)\Rightarrow z) - (x\Rightarrow(y\Rightarrow z))\nonumber\\
           &=&H(\mid x\mid,\mid y\mid,\mid z\mid)
           \sum\limits_{k=0}^{r}\sum\limits_{j=0}^{q} z_{1}\cdots z_{k}y_{1}\cdots y_{j}x y_{j+1} \cdots y_{q} z_{k+1}\cdots z_{r}\nonumber\\
           &&-H_2(\mid x\mid,\mid y\mid,\mid z\mid)
           \sum\limits_{k=0}^{r}\sum\limits_{j=0}^{k-1} z_{1}\cdots z_{j} x z_{j+1} \cdots z_{k}y z_{k+1} \cdots z_{r}\nonumber\\
           &&-H_2(\mid x\mid,\mid y\mid,\mid z\mid)
        \sum\limits_{k=0}^{r}\sum\limits_{j=k+1}^{r} z_{1}\cdots z_{k} y z_{k+1} \cdots z_{j}x z_{j+1} \cdots z_{r}\nonumber\\
           &:=&LHS,\label{2VDD-E2.3}
\end{eqnarray}
\begin{eqnarray}
           &&((y\Rightarrow x)\Rightarrow z) - (y\Rightarrow(x\Rightarrow z))\nonumber\\
           &=&H(\mid y\mid,\mid x\mid,\mid z\mid)
           \sum\limits_{k=0}^{r}\sum\limits_{j=0}^{p} z_{1}\cdots z_{k}x_{1}\cdots x_{j}y x_{j+1} \cdots x_{p} z_{k+1}\cdots z_{r}\nonumber\\
           &&-H_2(\mid y\mid,\mid x\mid,\mid z\mid)
           \sum\limits_{k=0}^{r}\sum\limits_{j=0}^{k-1} z_{1} \cdots z_{j}y z_{j+1} \cdots z_{k}x z_{k+1} \cdots z_{r}\nonumber\\
           &&-H_2(\mid y\mid,\mid x\mid,\mid z\mid)
           \sum\limits_{k=0}^{r}\sum\limits_{j=k+1}^{r} z_{1}\cdots z_{k} x z_{k+1} \cdots z_{j}y z_{j+1} \cdots z_{r}\nonumber\\
           &:=&RHS.\label{2VDD-E2.4}
\end{eqnarray}
  From (\ref{2VDD-E2.3}), we observe that the beginning term of $LHS$ is $yxz$ with the coefficient $H(\mid x\mid,\mid y\mid,\mid z\mid)$ and there is no other terms containing $yxz$. The remaining terms are those in which both $x$ and $y$ are inserted into $z$ and separated by at least one letter of $z$. \\
  Therefore, for $\forall$ $x$, $y$, $z$, $LHS= RHS$ if and only if  $H(\mid x\mid,\mid y\mid,\mid z\mid)= H(\mid y\mid,\mid x\mid,\mid z\mid)=0$ and $H_2(\mid x\mid,\mid y\mid,\mid z\mid)= H_2(\mid y\mid,\mid x\mid,\mid z\mid)$.
  We immediately obtain $(\ref{2VDD-E2.0})$.
\end{proof}
\begin{rem}
  If we denote $H_1(\mid x\mid,\mid y\mid,\mid z\mid)$ as $f(\mid x \mid,\mid y \mid)f(\mid x\mid+\mid y \mid,\mid z \mid)$, then we have
  $$H_1(\mid x\mid,\mid y\mid,\mid z\mid)=H_2(\mid x\mid,\mid y\mid,\mid z\mid),\ \ H_2(\mid x\mid,\mid y\mid,\mid z\mid)=H_2(\mid y\mid,\mid x\mid,\mid z\mid).$$ \\
  That is $$f(\mid x\mid, \mid y\mid) f(\mid x\mid+\mid y\mid,\mid z\mid)=f(\mid y \mid,\mid z\mid)f(\mid x \mid,\mid y\mid +\mid z \mid)$$ and
$$f(\mid y \mid,\mid z\mid)f(\mid x \mid,\mid y\mid +\mid z \mid)=f(\mid x \mid,\mid z\mid)f(\mid y \mid,\mid x\mid +\mid z \mid).$$
\end{rem}

\subsection{The existence of $f$}
  Recall that $f$ satisfies the following equations:
\begin{eqnarray}
f(m, n) f(m+n, p)= f(n, p) f(m, n+p),   \label{2VDD-E3.1}
\end{eqnarray}
\begin{eqnarray}
f(n,p) f(m, n+p)= f(m, p) f(n, m+p ),     \label{2VDD-E3.2}
\end{eqnarray}
where $ \forall \ m, n, p \in \mathbb{N}.$
\begin{example}
  We choose a symmetric bilinear function $g: \mathbb{N} \times \mathbb{N} \rightarrow \mathbb{R} $ given by $(m, n)\mapsto g(m, n)$. For example, $g(m, n)=k\cdot mn,$ where $k$ is a fixed positive number. Let $f(m, n)=\exp\{g(m, n)\}$. Then
\begin{eqnarray}
  &&f(m, n) f(m+n, p)\nonumber\\
  &=&\exp\{g(m, n)\} \exp\{g(m+n,p)\}\nonumber\\
  &=&\exp\{g(m, n)+ g(m+n,p)\}\nonumber\\
  &=&\exp\{g(m, n)+ g(m,p)+ g(n, p)\}\nonumber\\
  &=&\exp\{g(m, n+p)+ g(n, p)\}\nonumber\\
  &=&\exp\{g(m, n+p)\} \exp\{ g(n, p)\}\nonumber\\
  &=&f(m, n+p) f(n, p),\nonumber
\end{eqnarray}
\begin{eqnarray}
  &&f(n, p) f(m, n+p)\nonumber\\
  &=&\exp\{g(n, p)\} \exp\{g(m, n+p)\}\nonumber\\
  &=&\exp\{g(n, p)+ g(m,n+p)\}\nonumber\\
  &=&\exp\{g(m, n)+ g(m,p)+ g(n, p)\}\nonumber\\
  &=&\exp\{g(n, m+p)+ g(m, p)\}\nonumber\\
  &=&\exp\{g(n, m+p)\} \exp\{ g(m, p)\}\nonumber\\
  &=&f(n, m+p) f(m, p).\nonumber
\end{eqnarray}
Hence, $f$ satisfies (\ref{2VDD-E3.1}) and (\ref{2VDD-E3.2}).
\end{example}

\begin{example}
    Assume that the function $f: \mathbb{N}\times\mathbb{N}\longrightarrow \mathbb{R}$ is defined by
$f(m,n):=\frac{m!\,n!}{(m+n)!}$, where $m,n\in\mathbb{N}$. Then 
\begin{eqnarray}
    &&f(m,n)f(m+n,p)\nonumber \\
&=& \frac{m!\,n!}{(m+n)!}\cdot\frac{(m+n)!\,p!}{(m+n+p)!}\nonumber \\
&=& \frac{m!\,n!\,p!}{(m+n+p)!}\nonumber \\
&=& \frac{n!\,p!}{(n+p)!}\cdot\frac{m!\,(n+p)!}{(m+n+p)!}\nonumber \\
&=& f(n,p)f(m,n+p)\nonumber,
\end{eqnarray}
\begin{eqnarray}
    &&f(n, p) f(m, n+p)\nonumber\\
  &=&\frac{n!\,p!}{(n+p)!}\cdot\frac{m!\,(n+p)!}{(m+n+p)!}\nonumber \\
&=& \frac{m!\,n!\,p!}{(m+n+p)!}\nonumber \\
&=& \frac{m!\,p!}{(m+p)!}\cdot\frac{n!\,(m+p)!}{(m+n+p)!}\nonumber \\
&=& f(m,p)f(n,m+p).\nonumber
\end{eqnarray}
Hence, $f$ satisfies (\ref{2VDD-E3.1}) and (\ref{2VDD-E3.2}).
\end{example}

\begin{example}
\item[(a)] Let $f(m, n)\equiv C_0$, where $C_0$ is a constant. It is the same as the simplified insertion.
\item[(b)] We can easily check that $f$ satisfies the (\ref{2VDD-E3.1}) and (\ref{2VDD-E3.2}), where $f$ is defined as follows
\begin{eqnarray}f(m, n)=
\begin{cases}
1, &\text{if}\ m,n\ \mbox{are both odd}\cr 0, &\mbox{otherwise}\nonumber
\end{cases} 
\end{eqnarray}
\end{example}

\subsection{The symmetric of $f$}
   Let $n=p$ in (\ref{2VDD-E3.1}) and (\ref{2VDD-E3.2}). Then we get $f(m+p, p)= f(p, m+p),\ \mbox{for}\ m \in \mathbb{N}$. Hence, $$f(m, n)=f(n, m),\ \mbox{for}\ m, n \in \mathbb{N}.$$


\subsection{Application prospects in genetics}
The algebraic framework of left symmetric algebras constructed from DNA insertion operations provides insights into DNA recombination mechanisms. The operation of the modified insertion with functions $f(m,n)$ establishes an algebraic foundation for the quantitative study of length-dependent insertion, such as in mathematical descriptions of transposition processes. And the constraint condition $f(m, n) f(m+n, p)=f(n, p) f(m, n+p)=
f(m, p) f(n, m+p)$ provides a mathematical model for transposon expansion dynamics and selection pressure. When $f(m,n)$ takes the form of exponential decay, it naturally describes the phenomenon of fitness reduction caused by long-fragment genetic insertions. 
The result provides a mathematical basis for research areas in biological genetics and genomics in a certain way.
\bibliographystyle{plain}
\bibliography{ref}

\end{document}